\documentclass[10pt]{amsart}

\usepackage{amssymb,amscd}

\newtheorem{theorem}{Theorem}
\newtheorem*{claim}{Main~Claim}
\newtheorem{lemma}[theorem]{Lemma}

\numberwithin{equation}{section}
\numberwithin{theorem}{section}

\theoremstyle{definition}
\newtheorem{definition}[theorem]{Definition}

\newcommand{\A}{\mathcal{A}}
\newcommand{\B}{\mathcal{B}}
\newcommand{\C}{\mathcal{C}}
\newcommand{\D}{\mathcal{D}}
\newcommand{\e}{\varepsilon}
\newcommand{\E}{\mathcal{E}}
\newcommand{\h}{\mathfrak{H}}

\newcommand{\N}{\mathcal{N}}
\newcommand{\M}{\mathcal{M}}
\newcommand{\p}{\mathcal{P}}
\newcommand{\R}{\mathcal{R}}
\newcommand{\Z}{\mathcal{Z}}

\begin{document}
\title[Banach-Stone]{A new proof of the noncommutative Banach-Stone theorem}
\author{David Sherman}
\address{Department of Mathematics\\ University of California\\ Santa Barbara, CA 93106}
\email{dsherman@math.ucsb.edu}
\subjclass[2000]{Primary: 46B04; Secondary: 46L05}
\keywords{C*-algebra, isometry, Jordan isomorphism}

\begin{abstract}
Surjective isometries between unital C*-algebras were classified
in 1951 by Kadison \cite{K}.  In 1972 Paterson and Sinclair
\cite{PS} handled the nonunital case by assuming Kadison's theorem
and supplying some supplementary lemmas.  Here we
combine an observation of Paterson and Sinclair with variations on the methods of Yeadon \cite{Y} and the author \cite{S1}, producing a fundamentally new proof of the structure of surjective isometries between (nonunital) C*-algebras.  In the final section we indicate how our techniques may be applied to classify surjective isometries of noncommutative $L^p$ spaces, extending the main results of \cite{S1} to $0 < p \leq 1$.
\end{abstract}

\maketitle

\section{Introduction}

The main goal of this paper is to give a new proof of the
noncommutative Banach-Stone theorem, by which we mean a description of the surjective isometries between (possibly nonunital) C*-algebras.  This is accomplished in Sections 2 and 3; Section 4 indicates how the techniques of our proof can be used to produce new results about noncommutative $L^p$ spaces.  We begin with a little of the relevant history.

The first theorem of this type was proved by Banach and is as follows.

\begin{theorem} \cite[Theorem IX.4.3]{B}
Let $X$ and $Y$ be compact metric spaces, $C(X)$ and $C(Y)$ the
associated Banach spaces of real continuous functions equipped
with the sup norm.  If $T: C(X) \to C(Y)$ is a surjective
isometry, then there are a homeomorphism $\varphi: Y \to X$ and a
function $h \in C(Y)$ with $|h(y)| = 1$ for all $y \in Y$, such
that
\begin{equation} T(f) = h (f \circ
\varphi), \qquad \forall f \in C(X).
\end{equation}
\end{theorem}

Banach's result was improved by Stone \cite{St} to handle the case
where $X$ and $Y$ are compact Hausdorff spaces, and this is the
version commonly known as the ``Banach-Stone" theorem.  Later
authors extended the result to complex (and even vector-valued)
functions on more general spaces - see \cite[Chapters 1 and 2]{FJ}
for some of the details.

Abelian unital C*-algebras are exactly the (complex) algebras $C(X)$, $X$ compact and Hausdorff.  Removing the assumption of commutativity leads us to Kadison's 1951 noncommutative Banach-Stone theorem.

\begin{theorem} \label{T:k} \cite[Theorem 7]{K}
Let $T: \A \to \B$ be a surjective isometry between unital
C*-algebras.  Then there are a surjective Jordan *-isomorphism $J$
from $\A$ to $\B$ and a unitary $u \in \B$ such
that
\begin{equation} T(x) = uJ(x), \qquad \forall x \in \A.
\end{equation}
\end{theorem}

Other proofs of Theorem \ref{T:k} exist, and by now Banach-Stone-type theorems have been given for a full menagerie of algebraic structures: power algebras, Jordan algebras, Hilbert C*-modules... see \cite{W} and \cite[Chapter 6]{FJ} for examples and discussion of the literature.  (Perhaps the study of quantized Banach-Stone theorems is ``noncommutative geology"?)  As for \textit{nonunital} C*-algebras, the first satisfactory classification for surjective isometries was obtained in 1972 by Paterson and Sinclair
\cite{PS}.  (There is a related result in the 1969 Ph.D. dissertation of Harris \cite{H}.)  Paterson and Sinclair assumed Theorem \ref{T:k} and added a few elegant observations, producing

\begin{theorem} \label{T:ps} \cite[Theorem 1]{PS}
Let $T: \A \to \B$ be a surjective isometry between C*-algebras,
and let $M(\B)$ be the multiplier algebra of $\B$.  Then there are
a surjective Jordan *-isomorphism $J$ from $\A$ to $\B$ and a
unitary $u \in M(\B)$ such that
\begin{equation} T(x) = uJ(x),
\qquad \forall x \in \A.
\end{equation}
\end{theorem}

In this paper we give a new proof of Theorem \ref{T:ps}.  We do
not assume Theorem \ref{T:k} or make use of any of the techniques
involved in its proofs, but we do benefit from a lemma of
\cite{PS}.

\bigskip

Our proof proceeds in the following manner.  Given $T: \A \to \B$ a surjective isometry between C*-algebras, we consider the induced surjective isometries $T^*: \B^* \to \A^*$, $T^{**}: \A^{**} \to \B^{**}$. Since the
second dual of a C*-algebra is isometric to a von Neumann algebra, $T^*$ is a surjective isometry of preduals of von Neumann algebras.  We determine the
structure of $T^*$, and this reveals the structure of $T^{**}$.  Then a lemma from \cite{PS} allows us to describe the restriction of $T^{**}$ to $\A$, which is nothing but $T$.  We arrive precisely at Theorem \ref{T:ps}.

Thus the fundamental object here is $T^*$, and the focus of this note is a new way of deriving the structure of surjective isometries between von Neumann preduals.  If one has \textit{already} proved Theorem \ref{T:k}, then the predual result follows easily, since the dual of a surjective isometry between preduals is a surjective isometry between unital C*-algebras.  And without reliance on Theorem \ref{T:k}, even nonsurjective isometries between preduals have been described - see \cite[Lemma 3.6]{Ki} or \cite[Theorem 3.2]{S2}.  Our technique in this paper is unlike the two papers just mentioned, and there is some novelty in using the predual result to recover Theorem \ref{T:ps}.  In this way we birth a Banach-Stone theorem without explicitly using the geometry of the unit ball in a C*-algebra in order to pick out distinguished classes of operators.  (For example, Kadison characterized the extreme points \cite[Theorem 1]{K}.)  The replacement, at the level of duals, is a certain orthogonality relation.

Our derivation of the predual result can be modified to extend the main theorems of \cite{S1}.  There a description was obtained for surjective isometries of noncommutative $L^p$ spaces, $1 < p < \infty, \: p \neq 2$, but actually a significant part of that proof was originally done in a different manner.  The ``original'' version requires special considerations for finite type I summands, and was considered by the author to be less elegant overall, but it \textit{does} remain valid for $p \leq 1$.  In this paper we utilize the ideas of the original proof, applied concretely to $p=1$ and at least motivated for other values of $p$.

Let us be more explicit: for the benefit of readers who are primarily interested in the Banach-Stone result, we have put off all discussion of noncommutative $L^p$ spaces to Section \ref{S:Lp}.  Nonetheless the reader is apprised that some of the key steps in our proof of Theorem \ref{T:ps} are the case $p=1$ of known results about noncommutative $L^p$ spaces.  Since von Neumann preduals are relatively easy to work with - compared with general noncommutative $L^p$ spaces - we supply or sketch direct proofs.  Then in Section \ref{S:Lp} we indicate how the argument extends to $p \notin\{1,2\}$, allowing the main results of \cite{S1} to be generalized.  As much as possible we avoid repeating arguments from \cite{S1}.

For better or worse, this paper operates at two levels of sophistication.  On one hand we attempt to prove the Banach-Stone theorem with as few prerequisites as possible.  The required von Neumann algebra theory is classical, and we supply explanations for most of the nontrivial, post-1950s assertions which do not derive from \cite{S1}.  We would like to think that an operator algebraist from the 1950s could digest the proof.  But on the other hand we also have in mind a reader familiar with \cite{S1} (and the theory of noncommutative $L^p$ spaces).  The extension of the main theorems of \cite{S1} is a new result which has already had applications elsewhere \cite{HRR}. 

\section{Definitions, facts, lemmas} \label{S:facts}

\textbf{1.} For a C*-algebra $\A$, the \textit{multiplier algebra}
$M(\A)$ can be defined in several equivalent ways.  Abstractly, it
can be described as the largest C*-algebra in which $\A$ embeds as
an essential ideal.  For commutative C*-algebras, this is
equivalent to embedding the Gelfand spectrum in its Stone-\v{C}ech
compactification.

$M(\A)$ can also be constructed as the algebra of double
centralizers.  For a concrete realization, one takes a faithful
nondegenerate *-representation $\pi$ of $\A$ on a Hilbert space
$\h$. Then
$$M(\A) \simeq \{x \in \B(\h) \mid x\pi(\A) \subset \pi(\A), \: \pi(\A) x \subset \pi(\A) \} \subseteq \pi(\A)''.$$
A special case of this construction occurs when $\pi$ is the universal representation of $\A$, so that $\A \simeq \pi(\A) \subset \pi(\A)''$ can be isometrically identified by the 1954 Sherman-Takeda theorem \cite{Ta} with the canonical embedding $\A \hookrightarrow \A^{**}$.  With this identification, $M(\A)$ is the idealizer of $\A$ in $\A^{**}$.

Multiplier algebras of C*-algebras have been around since the 1960s and are discussed at length in
\cite[Chapter 2]{W-O}.  $M(\A)$ is always unital and equals $\A$ when $\A$ is itself unital.  As an example, $\B(\h)$ is the multiplier algebra of the C*-algebra of compact operators on $\h$. 

\bigskip

\textbf{2.} A linear map between operator algebras is \textit{Jordan}
when it preserves the Jordan product $(x, y) \mapsto (\frac{1}{2})(xy
+ yx)$.  This is equivalent to requiring that $J$ commute with squaring on self-adjoint (or all) elements, by distributivity and elementary algebra involving $(x + y)^2$.  A fundamental 1951 result of Kadison \cite[Theorem
10]{K} says that a surjective Jordan *-isomorphism between von Neumann algebras is the direct sum of a *-isomorphism and a *-antiisomorphism.

\bigskip

\textbf{3.} We will need one of the two lemmas with which Paterson and Sinclair paved their path from Theorem \ref{T:k} to Theorem \ref{T:ps}.  For the reader's convenience, we present a short proof (directly adapted from the original article).

\begin{lemma} \label{T:j} \cite[Lemma 2]{PS}
Let $K:\C \to \D$ be a Jordan *-monomorphism of C*-algebras, and assume $\D$ is unital.  If $v$ is
a unitary element of $\D$ such that $vK(\C)$ is a C*-subalgebra of $\D$, then $K(\C) = v K(\C)$.
\end{lemma}

\begin{proof} We will use that any C*-algebra $\E$ satisfies $\E^2 = \E$ (for example, by Cohen's 1959 factorization theorem \cite{C}).  First,
\begin{equation} \label{E:ps1}
vK(\C) = (v K(\C))(v K(\C)) = (K(\C)v^*)(v K(\C)) = (K(\C))^2.
\end{equation}
Now any element in $\C$ is a linear combination of four squares (because positive elements are squares).  Since $K$ commutes with squaring, \eqref{E:ps1} implies
\begin{equation} \label{E:ps2}
K(\C) \subseteq \text{span }(K(\C))^2 = \text{span } vK(\C) = vK(\C).
\end{equation}
We combine \eqref{E:ps1} and \eqref{E:ps2} to obtain the converse inclusion:
$$vK(\C) = (K(\C))^2 \subseteq (K(\C))(vK(\C)) = v^*(vK(\C))(vK(\C)) = v^*(vK(\C)) = K(\C).$$
\end{proof} 

\bigskip

\textbf{4.} Let $\M$ be a von Neumann algebra, and identify the
predual $\M_*$ with the normal linear functionals on $\M$.  Then $\M_*$ is an $\M-\M$ bimodule, with actions defined by
\begin{equation} \label{E:actions}
x\rho(\cdot) = \rho(\cdot x), \qquad \rho x(\cdot) = \rho(x \cdot),
\qquad x \in \M, \: \rho \in \M_*.
\end{equation}

\begin{lemma} \label{T:incl}
The two inclusions $\M \hookrightarrow \B(\M_*)$ as left or right multipliers are isometric, and their images are commutants of each other in $\B(\M_*)$.
\end{lemma}

\begin{proof} Clearly
\begin{equation} \label{E:cont}
\|x\rho y\| = \sup_{\|a\|=1} \rho(yax) \leq \|x\| \|\rho\| \|y\|, \qquad x,y \in \M, \: \rho \in \M_*,
\end{equation}
so the inclusions are norm-decreasing.  

To show the inclusion as right multipliers is isometric - the other inclusion being analogous - choose generic $x \in \M$ and $\e > 0$.  Let $v|x|$ be the polar decomposition of $x$, $p$ be the nonzero spectral projection of $|x|$ corresponding to $[\|x\| - \e, \|x\|]$, and $\varphi \in \M_*^+$ be a state with support $\leq p$.  Then
$$1 \geq \|\varphi v^*\| \geq \varphi v^*(v) \geq \varphi(p) = 1 \Rightarrow \|\varphi v^*\| = 1$$
and 
$$\|(\varphi v^*)x \| = \|\varphi |x| \| \geq \varphi(|x|) \geq \varphi((\|x\| - \e)p) = \|x\| - \e.$$
Since $\e$ is arbitrary, the operation of right multiplication by $x$ on $\M_*$ must have norm $\|x\|$.

For the latter assertion, assume $T \in \B(\M_*)$ commutes with the right action, and calculate
$$[T(\rho)](x) = [T(\rho)x](1) = [T(\rho x)](1) = [\rho x](T^*(1)) = [T^*(1) \rho](x)$$
for all $x \in \M$, $\rho \in \M_*$, so that $T$ is left multiplication by $T^*(1)$.  Of course the commutant of the left action can be calculated similary.
\end{proof}

The previous lemma is generalized in \cite[Lemma 1.1 and Corollary 1.6]{JS}, while the following lemma is a special case of \cite[Lemma 1.3]{JS}.  Recall that the strong topology is the point-norm topology that $\M$ acquires from any faithful normal *-representation on a Hilbert space.

\begin{lemma} \label{T:strong}
For a net $\{x_\alpha\}$ in the unit ball of a von Neumann algebra $\M$, the following are equivalent:
\begin{enumerate}
\item $x_\alpha \to 0$ strongly;
\item $x_\alpha \to 0$ in the point-norm topology that $\M$ acquires from its left (or right) action on $\M_*$;
\item (if $\M$ is $\sigma$-finite) $x_\alpha \psi \to 0$ in the norm topology of $\M_*$, for a single faithful $\psi \in \M_*^+$.
\end{enumerate}
\end{lemma}

\begin{proof} $(1) \Rightarrow (2)$: Any normal linear functional on $\M$ is a linear combination of four positive ones, so it suffices to consider an arbitrary state $\varphi \in \M_*^+$.  Let $\{\pi_\varphi, \h_\varphi, \xi_\varphi\}$ be the associated GNS representation.  Then
$$\|x_\alpha \varphi\| = \sup_{\|y\|=1} |\varphi(y x_\alpha)| = \sup_{\|y\| = 1} |<\pi_\varphi(x_\alpha) \xi_\varphi \mid \pi_\varphi(y^*) \xi_\varphi>| \leq \|\pi_\varphi(x_\alpha) \xi_\varphi\| \to 0.$$

$(2) \Rightarrow (1)$: Choose any $\xi \in \h$, where $\{\pi, \h\}$ is a normal *-representation of $\M$.  Let $\varphi_\xi \in \M_*^+$ be the associated vector functional, i.e. $\varphi_\xi(x) = <\pi(x) \xi \mid \xi>$.  Then
$$\|\pi(x_\alpha) \xi\|^2 = <\pi(x_\alpha^* x_\alpha) \xi \mid \xi> = \varphi_\xi(x_\alpha^* x_\alpha) = x_\alpha \varphi_\xi(x_\alpha^*) \leq \|x_\alpha \varphi_\xi\| \to 0.$$

$(2) \Rightarrow (3)$: Trivial.

$(3) \Rightarrow (2)$: We first claim that the subspace $\psi \M \subset \M_*$ is norm dense.  If not, Hahn-Banach guarantees the existence of a nonzero linear functional which vanishes on $\overline{\psi \M}$.  Such a functional arises from an element of $\M$, so we have a nonzero $x \in \M$ satisfying $\psi y (x) = 0$ for all $y \in \M$.  Letting $x = v|x|$ be the polar decomposition and setting $y = v^*$, we arrive at $\psi(|x|) = 0$, contradicting the faithfulness of $\psi$.

Now choose any $\rho \in \M_*$ and $\e > 0$.  Find $y \in \M$ with $\|\rho - \psi y\| < \e$.  Then 
$$\|x_\alpha \rho\| \leq \| x_\alpha \rho - x_\alpha \psi y \| + \| x_\alpha \psi y\| \leq \|\rho - \psi y \| + \|x_\alpha \psi \| \|y\|,$$
which is eventually less than $\e$.
\end{proof}

Each $\rho \in \M_*$ has a unique
polar decomposition as $v|\rho|$, where $|\rho| \in \M_*^+$ and $v$
is a partial isometry in $\M$ with $v^* v = s(|\rho|)$, the support
of $|\rho|$.  The left (resp. right) support of $\rho$
is denoted by $s_\ell(\rho)$ (resp. $s_r(\rho)$) and is equal to
$vv^*$ (resp. $v^*v$).  This goes back to 1950s work of Sakai, as do the earlier assertions in this subsection about preduals and independence of the strong topology.

We will need the case $p=1$ of the equality condition in the Clarkson inequality (see Theorem \ref{T:lpclarkson}), which we prove here using elementary techniques.  For positive functionals, Lemma \ref{T:clarkson} is shown in \cite[Theorem III.4.2(ii)]{T}.

\begin{lemma} \label{T:clarkson}
For a von Neumann algebra $\M$ and $\rho, \sigma \in \M_*$, we have
\begin{equation} \label{E:clarkson}
\|\rho + \sigma\| = \|\rho\| + \|\sigma\| = \|\rho - \sigma\| \iff
s_\ell(\rho) \perp s_\ell(\sigma), \quad s_r(\rho) \perp s_r(\sigma).
\end{equation}
\end{lemma}

\begin{proof}
Let $\rho = v|\rho|$ and $\sigma = w|\sigma|$ be the polar
decompositions.  If we assume the right-hand side of \eqref{E:clarkson}, then $v^*w = vw^* = 0$.  The left-hand side follows by evaluating $\rho \pm \sigma$ at $v^* \pm w^*$, which has norm 1.

Now assume the left-hand side of \eqref{E:clarkson}.  Necessarily we have $x,y$ in the unit ball of $\M$ with
$$(\rho + \sigma)(x) = \|\rho\| + \|\sigma\| = (\rho - \sigma)(y),$$
and this implies
$$\|\rho\| = \rho(x) = |\rho|(s(|\rho|) x v s(|\rho|)) \Rightarrow s(|\rho|) x v
s(|\rho|) = s(|\rho|) \Rightarrow s_r(\rho) x s_\ell(\rho) = v^*.$$
Similarly
$$s_r(\rho) y s_\ell(\rho) = v^*, \quad s_r(\sigma) x s_\ell(\sigma) = w^*, \quad s_r(\sigma) y s_\ell(\sigma) =
-w^*.$$

Write 
\begin{equation} \label{E:cl1} x = s_r(\rho) x s_\ell(\rho) + (1 - s_r(\rho))
x s_\ell(\rho) + s_r(\rho) x (1 - s_\ell(\rho)) + (1 - s_r(\rho)) x (1
- s_\ell(\rho)).
\end{equation}
Now $s_r(\rho) x s_\ell(\rho) = v^*$, so the inequalities
$$\|s_r(\rho)x x^* s_r(\rho)\| \leq 1, \qquad \|s_\ell(\rho) x^* x s_\ell(\rho)\| \leq 1,$$
force the two middle terms in \eqref{E:cl1} to drop out.  This leaves
\begin{equation} \label{E:cl2}
x = v^* + (1 - s_r(\rho)) x (1 - s_\ell(\rho)).
\end{equation}
Similarly
\begin{equation} \label{E:cl3}
y = v^* + (1
- s_r(\rho)) y (1 - s_\ell(\rho)).
\end{equation}
We may deduce from \eqref{E:cl2} and \eqref{E:cl3} that
\begin{equation} \label{E:cl4} s_r(\rho) x = v^* = s_r(\rho) y,
\end{equation}
and by an entirely analogous argument,
\begin{equation} \label{E:cl5}x s_\ell(\sigma) = w^* = - y s_\ell(\sigma).
\end{equation}

Finally \eqref{E:cl4} and \eqref{E:cl5} imply
$$s_r(\rho) x s_\ell(\sigma) = v^* s_\ell(\sigma) = s_r(\rho) y
s_\ell(\sigma) = - s_r(\rho) w^* = - s_r(\rho) x s_\ell(\sigma).$$
Then all terms in this equation are zero, and examination of the
second and fourth terms gives the conclusion.
\end{proof}

Following \cite{S1}, we say that functionals $\rho, \sigma$
satisfying \eqref{E:clarkson} are \textit{orthogonal}, and we write
$\rho \perp \sigma$.  From Lemma \ref{T:clarkson} it follows that
orthogonality is preserved by isometries between von Neumann
preduals.  Relations of this type have been exploited to study
isometries in a variety of contexts, starting with Banach \cite{B}
and developed especially by Lamperti \cite{L}.

We use the notion $\perp$ to define \textit{orthocomplements} as
well: for a set $S \subset \M_*$, $S^\perp$ is the set of elements
orthogonal to every element in $S$.

Finally, we use the notations $\p$ for ``projections of" and $\Z$ for ``center of", so for example $\p(\Z(\M))$ is the set of central projections in $\M$.  We also use $c(\cdot)$ for ``central support of", applied
to operators or elements in the predual.

\section{Proof of Theorem \ref{T:ps}} \label{S:proof}

Let $T: \A \to \B$ be a surjective isometry of C*-algebras.  For
brevity we set $\N = \A^{**}, \: \M = \B^{**}, \: \Phi = T^*$, so
that $\Phi:\M_* \to \N_*$ is a surjective isometry between von
Neumann preduals.

\begin{claim} There are a unitary $w \in \N$ and a surjective Jordan *-isomorphism $K: \M \to \N$ such that
\begin{equation} \label{E:claim}
\Phi(\rho) = w(\rho \circ K^{-1}), \qquad \rho \in \M_*.
\end{equation}
\end{claim}

Our proof of the Main Claim has three steps.  First we show that we may consider separately the finite type I summand of $\M$.  Then we prove the claim for finite type I algebras, starting with lines laid down by Yeadon \cite{Y}.  Finally we prove the claim for algebras with no finite type I summand, this time using a variation on the methods of \cite{S1}.

\bigskip

\textit{Step 1:}

\begin{lemma} \label{T:corner}
Let $\Phi$ be as above.  If $z \in
\mathcal{P}(\mathcal{Z}(\mathcal{M})),$ then
\begin{equation} \label{E:center}
\Phi(\M_* z) = \N_* z' \: \text{ for some } z' \in
\mathcal{P}(\mathcal{Z}(\mathcal{N})).
\end{equation}
The map $z \mapsto z'$ induces a surjective *-isomorphism from $\mathcal{Z}(\mathcal{M})$ to $\mathcal{Z}(\mathcal{N})$.
\end{lemma}

The proof is identical to Lemma 4.1 in \cite{S1}.  We continue to
use the apostrophe for the isomorphism map between centers.

Notice that Lemma \ref{T:corner} implies that $c(\Phi(\rho)) =
c(\rho)'$ for any $\rho \in \M_*$. Now for $z$ a central projection in an arbitrary von
Neumann algebra $\R$, define
$$N(z) = \sup \{n \mid \exists \rho_1, \rho_2, \dots \rho_n  \in \R_* \text{
with } c(\rho_j) = z, \: \rho_j \perp \rho_k \text{ for } j \neq
k\}.$$
Coming back to our context, we have that $N(z') = N(z)$.

If $z_k$ is the central projection onto the $\text{I}_k$ summand
in $\M$, $z_k$ is exactly characterized as the largest central
projection such that $N(z) = k, \: \forall z \leq z_k$.  It
follows that $z_k'$ is the central projection onto the
$\text{I}_k$ summand in $\N$.  Denoting the finite type I summand
of $\M$ as $\M_\text{I, fin}$, then, we have that $\Phi$ restricts
to a surjective isometry from $(\M_\text{I, fin})_*$ to
$(\N_\text{I, fin})_*$.

\bigskip

\textit{Step 2:} Yeadon \cite[Theorem 2]{Y} determined the form of all (not
necessarily surjective) isometries between noncommutative $L^p$
spaces ($1 \leq p < \infty, \: p \neq 2$) associated to semifinite
von Neumann algebras.  The first half of Step 2 is a variation of his method.

Let $\Phi:\M_* \to \N_*$ be a surjective isometry between preduals of finite type I von Neumann algebras.  Temporarily assume that $\M$ is $\sigma$-finite, and fix a faithful normal trace $\tau_\M$ on $\M$.  For each projection $p \in \p(\M)$, let $w_p \varphi_p$ be the polar decomposition of $\Phi(\tau_\M p)$.  We have that
$$w_1 \varphi_1 = \Phi(\tau_\M) = \Phi(\tau_\M p) + \Phi(\tau_\M (1-p)) = w_p \varphi_p + w_{(1-p)} \varphi_{(1-p)}.$$
Since $\tau_\M p$ and $\tau_\M (1-p)$ are orthogonal, so are their images.  Thus $w_p \varphi_{(1-p)} = w_{(1-p)} \varphi_p = 0$, and
$$w_1 \varphi_1 = (w_p + w_{(1-p)}) (\varphi_p + \varphi_{(1-p)}).$$
Moreover $w_p + w_{(1-p)}$ is still a partial isometry, as the summands have orthogonal left and right supports, so uniqueness of the polar decomposition implies
\begin{equation} \label{E:yea}
w_p + w_{(1-p)} = w_1, \qquad \varphi_1 = \varphi_p + \varphi_{(1-p)}.
\end{equation}  

Define
$$K:\p(\M) \to \p(\N); \qquad p \mapsto s_r(\Phi(\tau_\M p)) = s(\varphi_p).$$
We have by \eqref{E:yea} that
\begin{equation} \label{E:yea2}
K(p) \varphi_1 = \varphi_p = \varphi_1 K(p) \quad \text{and} \quad \Phi(\tau_\M p) = w_1 \varphi_1 K(p), \qquad p \in \p(\M).
\end{equation}
We already know that $K$ is additive on orthogonal projections; extend $K$ to real linear combinations of orthogonal projections, then by continuity to arbitrary self-adjoint operators, then by complex linearity to all of $\M$.  Because of \eqref{E:yea2} we have
\begin{equation} \label{E:commute}
K(x) \varphi_1 = \varphi_1 K(x) \quad \text{and} \quad \Phi(\tau_\M x) = w_1 \varphi_1 K(x), \qquad x \in \M,
\end{equation}
which guarantees that $K$ is a well-defined injective linear map.  By construction $K$ is *-preserving and commutes with squaring on self-adjoint elements, so it is also Jordan.

The density of $\tau_\M \M$ in $\M_*$ implies the density of $w_1 \varphi_1 K(\M)$ in $\N_*$.  Since every element of $w_1 \varphi_1 K(\M)$ vanishes on $(1 - s_\ell(w_1))$, we must have $s_\ell(w_1) = 1$.  The finiteness of  $\N$ then means that $w_1$ is unitary and $\varphi_1$ is faithful.

For a bounded net $\{x_\alpha\} \subset \M$, the convergence of $\tau_\M x_\alpha$ is equivalent to the convergence of $\varphi_1 K(x_\alpha)$.  Lemma \ref{T:strong} tells us that $K$ is strongly continuous on the unit ball of $\M$, so that the unit ball of $K(\M)$ is strongly and weakly closed in $\N$.  By the Krein-Smullyan theorem $K(\M)$ is weakly closed in $\N$. 

We claim that $K(\M)$ is also weakly dense in $\N$ and so must equal $\N$.  It is enough to show that each element of $\N_*$ attains its norm when restricted to $K(\M)$, and by density of $\varphi_1 K(\M)$ in $\N_*$ it suffices to prove this for a functional of the form $\varphi_1 K(x)$, $x \in \M$.  In fact the norm is attained at $K(v^*)$, where $x = v |x|$ is the polar decomposition:
\begin{align*}
[\varphi_1 K(x)] (K(v^*)) &= \varphi_1(K(x) K(v^*)) \overset{\eqref{E:commute}}{=} \varphi_1 \left( \frac{K(x)K(v^*)+K(v^*)K(x)}{2} \right) \\ &= \varphi_1 \left(K \left( \frac{x v^* + v^* x}{2} \right) \right) = \varphi_1 \left(K \left( \frac{|x^*| + |x|}{2}\right) \right) \\ &= \left\|\varphi_1 K\left( \frac{|x^*| + |x|}{2} \right) \right\|_{\N_*} = \left\| w_1 \varphi_1 K\left( \frac{|x^*| + |x|}{2} \right) \right\|_{\N_*}\\ &= \left\|\tau_\M \left(\frac{|x^*| + |x|}{2}\right) \right\|_{\M_*} = \| \tau_\M x \|_{\M_*} = \|\varphi_1 K(x)\|_{\N_*}.
\end{align*}
This establishes that $K$ is surjective.

By \eqref{E:commute}, $\varphi_1$ is therefore a finite faithful trace on $\N$.  Actually $\varphi_1 =\tau_\M \circ K^{-1}$, since for any $h \in \N_+$,
$$\tau_\M (K^{-1}(h)) = \|\tau_\M K^{-1}(h) \|_{\M_*} = \| \Phi(\tau_\M K^{-1}(h))\|_{\N_*} = \|w_1 \varphi_1 h \|_{\N_*} = \varphi_1(h).$$
We have therefore that
$$\Phi(\tau_\M x) = w_1 (\tau_\M \circ K^{-1}) K(x) = w_1 ((\tau_\M x) \circ K^{-1}), \qquad x \in \M.$$
(The second equality uses again that the Jordan product of two operators has the same trace as the usual product.)  By density of $\tau_\M \M$ we may conclude \eqref{E:claim} for all $\rho \in \M_*$, taking $w = w_1$.

If $\M$ is not $\sigma$-finite, we apply the above argument to each of its $\sigma$-finite central summands.  Since the isometries agree on their intersections, so do the associated Jordan *-isomorphisms and unitaries (partial isometries in $\M$).  It follows that there are a global Jordan  *-isomorphism and unitary of which these are restrictions, and \eqref{E:claim} holds in general.

\bigskip

\textit{Step 3:} Now we assume that $\Phi: \M_* \to \N_*$ is a surjective isometry between preduals of algebras which have no finite type I summand.  We start with a useful

\begin{definition}  \cite{S1} Let $\M$ be a von Neumann algebra.  A subspace of $\M_*$ is called a \textit{corner} if it is of the
form $q_1 \M_* q_2$ for some $q_1,q_2 \in \p(\M)$.  Corners with
$q_1 = 1$ (resp. $q_2 = 1$) are called \textit{columns} (resp.
\textit{rows}).  Notice that a corner has a unique representation
in which $c(q_1) = c(q_2)$.  A corner of the form $\M_* z$, $z \in \p(\Z(\M))$, is called a \textit{central summand} (as is the algebra $\M z$).
\end{definition}

Statements (3)-(5) of the following lemma are also included in \cite[Lemma 3.1]{S1}.

\begin{lemma} \label{T:perp} ${}$
Let $\M$ and $\N$ be von Neumann algebras.
\begin{enumerate}
\item When $q_1, q_2 \in \p(\M)$ satisfy $c(q_1) = c(q_2)$, then
$(q_1 \M_* q_2)^\perp = (1-q_1)\M_*(1-q_2)$. \item A corner in
$\M_*$ has orthocomplement $\{0\}$ if and only if it can be
(re)written as $\M_* r_1 z + r_2 \M_* (1-z)$, where $z \in
\p(\Z(\M))$, $r_1,r_2 \in \p(\M)$, and $c(r_1 z + r_2 (1-z)) = 1$.
\item If $\Phi: \M_* \to \N_*$ is a surjective isometry and $S
\subset \M_*$, then $T(S^\perp) = T(S)^\perp$. \item The
intersection of any collection of corners in $\M_*$ is a corner.
\item For any set $S \subset \M_*$, $S^\perp$ is a corner. \item
The closure of the union of an increasing net of corners in $\M_*$
is a corner.
\end{enumerate}
\end{lemma}

\begin{proof}
The first statement is obvious and implies the second.  By \eqref{E:clarkson}, $\Phi$ and
$\Phi^{-1}$ preserve orthogonality, proving the third statement.
For the fourth, let $\{p_\alpha\}, \{q_\alpha\} \subset \p(\M)$;
then
$$\bigcap p_\alpha \M_* q_\alpha = (\wedge
p_\alpha) \M_* (\wedge q_\alpha).$$ The fifth follows from noting
that $\{\rho\}^\perp = (1-s_\ell(\rho)) \M_* (1-s_r(\rho))$ and
applying the fourth to the expression
$$S^\perp = \bigcap_{\rho \in S} \{\rho\}^\perp.$$

To prove the sixth, assume that $\{p_\alpha\}, \{q_\alpha\}$ are
increasing nets in $\p(\M)$.  Necessarily we have the strong
convergences $p_\alpha \to p = (\sup p_\alpha)$, $q_\alpha \to q =
(\sup q_\alpha)$.  Then for any $\rho \in \M_*$, \eqref{E:cont} and Lemma \ref{T:strong} imply
\begin{align*}
\|p \rho q -
p_\alpha \rho q_\alpha\| \leq \|p \rho (q - q_\alpha) \| + \|[(p -
p_\alpha) \rho] q_\alpha \| \leq \| \rho (q - q_\alpha) \| + \|(p -
p_\alpha) \rho \| \to 0.
\end{align*}
It follows that
$$\overline{\bigcup p_\alpha \M_* q_\alpha} = p
\M_* q. \qedhere$$

\end{proof}

A version of the next lemma was proved in
\cite[Lemma 4.3]{S1} for noncommutative $L^p$ spaces, $1 < p < \infty$ only, using different techniques.

\begin{lemma} \label{T:allcorner} Assume $\M$ and $\N$ have no finite
type I summand, and let $\Phi: \M_* \to \N_*$ be a
surjective isometry.
\begin{enumerate}
\item $\Phi$ takes corners to corners. \item If $q \in
\mathcal{P}(\mathcal{M})$ satisfies $c(q) = c(1-q) = 1$, then
$$\Phi(\mathcal{M}_*q) = \mathcal{N}_* q_1 z' + q_2
\mathcal{N}_*(1-z'),$$
for some $q_1,q_2 \in \mathcal{P}(\mathcal{N}),$ $z'
\in \mathcal{P}(\mathcal{Z}(\mathcal{N})),$ with $q_1 z' + q_2 (1-z')$ strictly between 0 and 1 on every central summand.  (More technically, $c(q_1 z' + q_2
(1-z')) = c(1 - (q_1 z' + q_2 (1-z')) = 1$.)
\end{enumerate}
\end{lemma}

\begin{proof}
Let $p_1 \M_* p_2$ be a corner with $c(p_1) = c(p_2)$, and first
assume that $c(1- p_1) = c(1- p_2)$.  In this case Lemma
\ref{T:perp}(1,3,5) tells us that $p_1 \M_* p_2$ and $(1 - p_1)
\M_* (1 - p_2)$ are orthocomplements of each other, their
images are orthocomplements as well, and therefore the images are
corners.

If $c(p_1) = c(p_2)$ but $c(1- p_1) \neq c(1- p_2)$, then the
corner contains a nonzero column or row.  Since there is no
finite type I summand, such a corner can be written as the closure
of an increasing union of corners covered by the first paragraph.
By Lemma \ref{T:perp}(6), the image is a corner.

The second statement is a consequence of the first part and Lemmas
\ref{T:perp}(2) and \ref{T:corner}.  Since $\M_* q$ contains no central summand and $(\M_* q)^\perp = \{0\}$, the
same holds for its image, whence $q_1 z' + q_2 (1-z')$ is strictly
between 0 and 1 on every central summand.
\end{proof}

At this point we can apply the same arguments as those given in Section 4 of \cite{S1}, from Lemma 4.4 until the end.  These arguments are self-contained, except for references to \cite{JS} which are covered in our context by Lemma \ref{T:incl} of this paper, and one assertion mentioned in the second item below.  The main points are these:

\begin{itemize}
\item Any choice of $q$ in Lemma \ref{T:allcorner}(2) produces the same central projection $z'$.
\item On $\M_* z$, $\Phi$ takes columns to columns; this induces an orthogonality-preserving map $\pi$ between projection lattices via $\Phi(\M_* z q) = \N_* z' \pi(q)$.  The map $\pi$ extends to a *-isomorphism between $\M z$ and $\N z'$ which is also an intertwiner:
$$\Phi(\rho x) = \Phi(\rho) \pi(x), \qquad \rho \in \M_*z, \: x \in \M z.$$
(In the context of \cite{S1} the fact that $\pi$ preserves orthogonality of projections is justified by properties of the semi-inner product.  More generally one can use that $T$ preserves orthogonality of predual/$L^p$ vectors to obtain that $\pi$ preserves orthogonality of $\sigma$-finite projections, then pass to general projections by considering increasing nets.)
\item On $\M_* z$, the map $\Phi$ can be decomposed as
$$\Phi(\rho) = w_1(\rho \circ \pi^{-1}), \qquad \rho \in \M_*z,$$
where $w_1$ is a unitary in $\M z$.
\item On $\M_* (1-z)$, $\Phi$ takes columns to rows, the analogous map $\bar{\pi}$ extends to a *-antiisomorphism, and 
$$\Phi(\rho) = w_2(\rho \circ \bar{\pi}^{-1}), \qquad \rho \in \M_*(1-z),$$
for a unitary $w_2$ in $\M(1-z)$.
\item Taking $w = w_1 + w_2$, $K = \pi \oplus \bar{\pi}$, we obtain \eqref{E:claim}.
\end{itemize}

This ends Step 3.  The Main Claim is therefore established by considering the restrictions of $\Phi$ to the finite type I summand of $\M$ and its complement, then adding the Jordan isomorphisms and unitaries.  (The unitaries in the two summands add to a unitary in $\M$.)

\bigskip

Having proved the Main Claim, decompose the surjective Jordan *-isomorphism $K$ as the sum of the *-isomorphism $K_1$ on $\M z$ and the *-antiisomorphism $K_2$ on $\M (1-z)$.  We have, for any $\rho \in \M_*$, $y \in \N$,
\begin{align} \label{E:form} \Phi(&\rho)(y) = [w(\rho \circ K^{-1})](y) \\ \notag &=  \rho(K_1^{-1}(y w z') + K_2^{-1}(yw (1-z'))) \\ \notag &= \rho((K_1^{-1}(wz') + K_2^{-1}(w(1-z')))(K_1^{-1}(w^* y w z') + K_2^{-1}(y(1-z'))) \\ \notag &= \rho(u J(y)).
\end{align}
Here $u$ is the unitary $(K_1^{-1}(wz') + K_2^{-1}(w(1-z')))$ in $\M$, and $J: \M \to \N$ is the surjective Jordan *-isomorphism $K^{-1} \circ \text{Ad }(wz' + (1-z'))$.

Now we return to the original map $T:\A \to \B$.  From \eqref{E:form}
we have that 
$$\rho(\Phi^*(y)) = [\Phi(\rho)](y) = [(\rho u) \circ J](y) = \rho (u J(y)), \qquad \rho \in \M_*, \: y \in \N,$$
so that $T^{**}(y) = \Phi^*(y) = u J(y)$.  Of course $T = T^{**} \mid_\A$.  

We apply Lemma \ref{T:j}, taking $\C = \A$, $\D = \B^{**}$, $v = u$, and $K = J\mid_\A$.  Since $uJ(\A) = T^{**}(\A) = T(\A) = \B$ is a C*-subalgebra of $\B^{**}$, we conclude by the lemma that $J(\A) = uJ(\A)$, whence $u \B = u(uJ(A)) = uJ(A) = \B$.  We also have that $\B u = [u^* \B]^* = [u^* (u \B)]^* = \B$, so that $u$ is a multiplier of $\B$.  The proof of
Theorem \ref{T:ps} is complete.

\section{Extension to noncommutative $L^p$ spaces} \label{S:Lp}

The main vehicle in our proof of Theorem \ref{T:ps} is the description of surjective isometries of preduals given in the Main Claim.  As mentioned in the introduction, the description itself is not new.  But our method of proof for preduals (noncommutative $L^1$ spaces) can be adapted to describe surjective isometries of noncommutative $L^p$ spaces, $0 < p < \infty, \, p \neq 2$, and this \textit{is} new when $p < 1$.  In this section we briefly describe the necessary changes and corresponding results.  We do not define noncommutative $L^p$ spaces formally; the reader seeking background may wish to consult the recent exposition \cite{PX}.  

In Subsection 2.4, all statements about $\M_*$ have obvious translations in terms of $L^p(\M)$, except that the bimodule structure is not as simple as \eqref{E:actions}, and Lemma \ref{T:clarkson} becomes
\begin{theorem}  (Equality condition for noncommutative Clarkson inequality) \label{T:lpclarkson}

For $\xi, \eta \in L^p(\mathcal{M})$, $0 < p < \infty$, $p \ne 2$,
\begin{equation}
\|\xi + \eta\|^p + \|\xi - \eta\|^p = 2(\| \xi \|^p + \| \eta
\|^p) \iff \xi \eta^* = \xi^* \eta = 0.
\end{equation}
\end{theorem}
\noindent Theorem \ref{T:lpclarkson} is due to Raynaud and Xu \cite{RX}, generalizing an earlier theorem of Kosaki \cite{Ko}.  Because of it, the orthogonality relation is preserved by isometries, and the ``calculus" of corners does not depend on the value of $p$.

The proof in Section \ref{S:proof} also proceeds as before, changing $\M_*$ to $L^p(\M)$.  However, \eqref{E:claim} can be stated for \textit{positive} $L^p$ elements only.  (And this is enough to describe a linear map, as any $L^p$ element is a linear combination of four positive ones.)  Positive $L^p$ elements may be viewed as $p$th roots of positive $L^1$ elements, and using this as a basis for notation, \eqref{E:claim} becomes
\begin{equation}
\Phi(\varphi^{1/p}) = w(\varphi \circ K^{-1})^{1/p}, \qquad \varphi \in \M_*^+.
\end{equation}

The Main Claim extends \cite[Theorem 1.2]{S1}, which only covered $p>1$, as follows:

\begin{theorem} (Noncommutative $L^p$ Banach-Stone theorem)\\
Let $T:L^p(\mathcal{M}) \to L^p(\mathcal{N})$ be a surjective isometry, where $\mathcal{M}$ and $\mathcal{N}$ are von Neumann algebras and $0 < p < \infty$, $p \ne 2$.  Then there are a surjective Jordan *-isomorphism $J:\mathcal{M} \to \mathcal{N}$ and a unitary $u \in \mathcal{N}$ such that
\begin{equation}
T(\varphi^{1/p}) = u(\varphi \circ J^{-1})^{1/p}, \qquad \forall \varphi \in \mathcal{M}_*^+.
\end{equation}
\end{theorem}

It is also possible to deduce the corresponding extension of \cite[Theorem 1.1]{S1}, which was previously stated for $p \geq 1$:

\begin{theorem} \label{T:s1} Let $\mathcal{M}$ and $\mathcal{N}$ be von Neumann algebras, and $0 < p \le \infty$, $p \ne 2$.  The following are equivalent:
\begin{enumerate}
\item $\mathcal{M}$ and $\mathcal{N}$ are Jordan *-isomorphic;
\item $L^p(\mathcal{M})$ and $L^p(\mathcal{N})$ are isometrically isomorphic as Banach spaces (or $p$-Banach spaces, when $p<1$).
\end{enumerate}
\end{theorem}

\end{document}